\documentclass[12pt]{article}

\usepackage{amssymb,amsmath,amsbsy, amsthm}
\usepackage{amscd}
\usepackage{pb-diagram}
\usepackage{pstricks,pst-node}
\usepackage{amsfonts}
\usepackage{latexsym}
\usepackage{pstricks,pst-node}
\usepackage{tikz}
\usepackage{fullpage}
\usepackage{xcolor}

\usepackage{hyperref}
\hypersetup{
     colorlinks   = true,
     citecolor    = blue,
    linkcolor = blue
}

\swapnumbers\newtheorem{theorem}{Theorem}[section]
\newtheorem{lemma}[theorem]{Lemma}
\newtheorem{corollary}[theorem]{Corollary}

\newtheorem{proposition}[theorem]{Proposition}
\newtheorem{example}[theorem]{Example}

\theoremstyle{definition}

\newtheorem{remark}[theorem]{Remark}

\numberwithin{equation}{section}



\newcommand{\diam}{\mathrm{diam}}

\def\RR{{\hbox{{\rm I$\!$\rm R}}}}

\def\NN{{\hbox{{\rm I$\!$\rm N}}}}

\def\tto{\ \hbox{$\to\!\!\!\!\!\to$}\ }

\newcommand{\A}{\mathcal{A}}

\newcommand{\K}{\mathcal{K}}

\newcommand{\M}{\mathcal{M}}
\newcommand{\T}{\mathcal{T}}
\newcommand{\U}{\mathcal{U}}
\newcommand{\C}{\mathcal{C}}
\newcommand{\D}{\mathcal{D}}

\newcommand{\FF}{\mathcal{F}}
\newcommand{\G}{\mathcal{G}}
\newcommand{\HH}{\mathcal{H}}

\newcommand{\V}{\mathcal{V}}

\newcommand{\Q}{\rightline{\bf Q.E.D.}\vskip.1in}

\begin{document}

\date{}
\title{More on $\T$-closed sets}

\author{Javier Camargo and Sergio Mac\'{\i}as
\footnote{{\it 2020 Mathematics Subject Classification}.
Primary: 54B20, 54C60, 54F15.
\newline
\textit{Key words and phrases.} Composant, continuum, continuum of 
colocal connectedness, continuumwise connected space,
hyperspace, property of Kelley, not a strong centre, pseudo-arc,
set function $\T$, shore subcontinuum, 
strongly continuumwise connected space, $\T$-closed set,
$\T$-closed subcontinuum.\hfil\break
}}

\maketitle

\begin{abstract}
We consider properties of the diagonal of a continuum that are used 
later in the paper. We continue the study of $\T$-closed subsets of a 
continuum $X$. We prove that
for a continuum $X$, the statements: $\Delta_X$ is a nonblock subcontinuum of $X^2$,
$\Delta_X$ is a shore subcontinuum of $X^2$ and $\Delta_X$ is not a strong 
centre of $X^2$ are equivalent, this result
answers in the negative Questions 35 and 36 and 
Question 38 ($i\in\{4,5\}$) of the paper ``{\sl Diagonals on the edge of the 
square of a continuum}, by A. Illanes, V. Mart\'{\i}nez-de-la-Vega, 
J. M. Mart\'{\i}nez-Montejano and D. Michalik''.  We also include
an example,
giving a negative answer to Question 1.2 of the paper
``{\sl Concerning when $\FF_1(X)$ is a continuum of colocal connectedness in 
hyperspaces and symmetric products}, Colloquium Math., 160 (2020), 297-307'',
by V. Mart\'{\i}nez-de-la-Vega, J. M. Mart\'{\i}nez-Montejano.
We characterised the $\T$-closed subcontinua of the square of the
pseudo-arc. We prove that the $\T$-closed sets of the product of two continua is compact
if and only if such product is locally connected. We show that 
for a chainable continuum $X$, $\Delta_X$ is a $\T$-closed subcontinuum of
$X^2$ if and only if $X$ is an arc. We prove that if $X$ is a continuum with the property of Kelley, then the
following are equivalent: $\Delta_X$ is a $\T$-closed subcontinuum of $X^2$,
$X^2\setminus\Delta_X$ is strongly continuumwise connected, $\Delta_X$ is a
subcontinuum of colocal connectedness, and $X^2\setminus\Delta_X$ is
continuumwise connected. We give models for the families of $\T$-closed sets and $\T$-closed subcontinua
of various families of continua.
\end{abstract}

\section{Introduction}\label{Intro}

Professor F. Burton Jones defined the
set functions $\T$ in \cite{J3}. Since then many 
properties related with this function have been studied.

The first use of the set functions $\T^n$ to study continua began in 1962 and 
may be found in \cite{DSS1} and \cite{DSS2}. Later, in 1967, the set function
$\T$ was applied to give sufficient conditions for the invertibility of continua
\cite{DD}. In 1968, this set function was used to study connectedness im kleinen
of continua \cite{Davis1}. In 1970, properties of continua for which the set 
function $\T$ is continuous \cite{B1}. In
1971, continuum neighborhoods and
filterbases were used to generalize the additivity of the set function $\T$ \cite{BD}.
In 2017, a long standing question was solved, it was shown that the continuity of 
the set function $\T$ implies the additivity of $\T$ \cite{CU}. In 2021, a book gathering 
many of the properties of the set function $\T$ was published \cite{MS}.

We begin considering properties of the diagonal of a continuum that are used 
later in the paper. Then we obtain properties of $\T$-closed sets in a general
continuum and we see when the diagonal of the square of a continuum is or is not
a $\T$-closed subcontinuum of the square of that continuum.
The paper is divided in six sections. After this introduction, in section~\ref{Def}, we give the
appropriate definitions of the concepts need to read the paper. In particular, if $X$ is a 
continuum and $A$ is a subcontinuum of $X$, then the notions of $A$ is a subcontinuum
of colocal connectedness, $X\setminus A$ is strongly continuumwise connected,
$X\setminus A$ is continuumwise connected, $A$ is a nonblock subcontinuum,
$A$ is a shore subcontinuum and $A$ is not a strong centre are introduced. Also,
some preliminary results related to these concepts are proved. In section~\ref{Diagonal},
we present properties of the diagonal of the square of a continuum $X$ related to
the concepts introduced in the previous section. For example, we prove that if $X$ is
a $\theta_1$-continuum, then the diagonal, $\Delta_X$, of $X^2$ is not a strong
centre of $X^2$ (Theorem~\ref{theta1cont}). Also, we show that if $X$ is a
continuum containing a terminal chainable subcontinuum, then $X^2\setminus\Delta_X$
is not continuumwise connected (Theorem~\ref{chaintermsubcont}). We prove that
for a continuum $X$, the statements: $\Delta_X$ is a nonblock subcontinuum of $X^2$,
$\Delta_X$ is a shore subcontinuum of $X^2$ and $\Delta_X$ is not a strong 
centre of $X^2$ are equivalent (Corollary~\ref{scimpliesnonblock2}), this 
result answers in the negative \cite[Questions 35 and 36]{IMMM2} and 
\cite[Question 38]{IMMM2}, for $i\in\{4,5\}$. We also present an 
example providing a negative answer to \cite[Question 1.2]{MM} 
(Proposition~\ref{examnegMM}). In section~\ref{Tclosed}, we give
new properties of $\T$-closed sets. We prove a characterization of $\T$-closed
sets (Theorem~\ref{charTclosedset}). We observe that  \cite[Problem 6.7]{CCOR} and 
\cite[Question 3.23]{OPQV} have partial answers when the continuum has the
property of Kelley and a positive answer when the continuum is locally connected
(Remarks~\ref{rem001} and~\ref{rem002}, respectively). We characterised the
family of $\T$-closed subcontinua for the square of the pseudo-arc
(Theorem~\ref{theo0003}). We prove that the family of $\T$-closed sets
($\T$-closed subcontinua) of the product of two continua is compact if and
only if such product is locally connected (Theorem~\ref{TclosedcompactiffXYlocconn}
and Corollary~\ref{TclosedcompactiffXYlocconn2}, respectively). We also
show that the family of $\T$-closed subcontinua of the pseudo-arc is
arcwise connected but not compact (Theorem~\ref{pseudoarcTclosedconn}).
We prove that if $X$ is a continuum and $A$ is a $\T$-closed subcontinuum of $X$, 
then the notions of $A$ is a subcontinuum
of colocal connectedness, $X\setminus A$ is strongly continuumwise connected,
$X\setminus A$ is continuumwise connected, $A$ is a nonblock subcontinuum,
$A$ is a shore subcontinuum and $A$ is not a strong centre are equivalent
(Theorem~\ref{ATclosedequiv}). In section~\ref{DeltaTclosed}, we restrict our study
to properties when the diagonal of the square of a continuum is or is not $\T$-closed.
We show that if $X$ is a chainable continuum containing a terminal subcontinuum,
then $\Delta_X$ is not a $\T$-closed subcontinuum of $X^2$ 
(Theorem~\ref{ctheochainterminal}). We prove that for a chainable continuum $X$,
$\T(\Delta_X)=X^2$ if and only if $X$ is indecomposable (Corollary~\ref{XchainableDeltaX2}).
Also, for a chainable continuum $X$, $\Delta_X$ is a $\T$-closed subcontinuum of
$X^2$ if and only if $X$ is an arc (Theorem~\ref{XchainDeltaTclosedXarc}).
We show that if $X$ is a continuum with the property of Kelley, then the
following are equivalent: $\Delta_X$ is a $\T$-closed subcontinuum of $X^2$,
$X^2\setminus\Delta_X$ is strongly continuumwise connected, $\Delta_X$ is a
subcontinuum of colocal connectedness, and $X^2\setminus\Delta_X$ is
continuumwise connected (Theorem~\ref{Deltachar}). In section~\ref{Models},
we give models for the families of $\T$-closed sets and $\T$-closed subcontinua
for continuously irreducible continua, continuously
type $A'$ $\theta$-continua, and continuously type $A'$ $\theta_n$-continua
(Theorem~\ref{models}). Also, we present models
for such families for the Menger universal curve of pseudo-arcs, the Sierpi\'nski
universal plane curve of pseudo-arcs, the circle of pseudo-arcs and the arc of
pseudo-arcs (Theorem~\ref{modelshomog}).

\section{Definitions}\label{Def}

Given a metric space $Z$ and a nonempty subset $A$ of $Z$,
$Int(A)$ and $Cl(A)$ denote the interior and the closure of $A$ in
$Z$, respectively. If $A$ is a nonempty closed subset of $Z$
and $\varepsilon>0$, then 
$\V_\varepsilon(A)=\{z\in Z\ |\ d(z,A)<\varepsilon\}$.

A continuous function between compact metric spaces is a {\it map}.
A surjective map 
$f\colon X\tto Y$ between compact metric spaces is said to be:
\medskip

$\bullet$ {\it atomic} provided that for each subcontinuum $L$ of $X$ such
that $f(L)$ is nondegenerate, we have that $L=f^{-1}(f(L))$.
\smallskip

$\bullet$ {\it closed} if $f(C)$ is closed in $Y$ for every 
closed subset $C$ of $X$.
\smallskip

$\bullet$ {\it confluent} provided that for each subcontinuum $Q$ of $Y$,
and every component $K$ of $f^{-1}(Q)$, we have that $f(K)=Q$.
\smallskip

$\bullet$ {\it monotone} if $f^{-1}(B)$ is connected for each connected
subset $B$ of $Y$;
\smallskip

$\bullet$ {\it open} provided that $f(U)$ is open in $Y$ for every open
subset $U$ of $X$.
\medskip

Given a metric space $Z$, a {\it decomposition} of $Z$ is a family $\G$ of nonempty
and mutually disjoint subsets of $Z$ such that $\bigcup\G=Z$. A decomposition
$\G$ of a metric space $Z$ is said to be {\it upper semicontinuous} if the quotient map
$q\colon Z\tto Z/\G$ is closed. The decomposition is {\it continuous} provided that
the quotient map is both closed and open.

A {\it continuum} is a compact, connected, metric space. A
{\it subcontinuum} is a continuum contained in a metric space.
Let $X$ be a continuum and let $p$ be an element of $X$.
Then the {\it composant of $p$ in $X$} is
$$\kappa(p)=\bigcup\{K\ |\ K\ \hbox{is a subcontinuum of $X$ and}\ p\in K\}.$$
A continuum $X$ is {\it decomposable} provided that there exist
two proper subcontinua $K$ and $L$ such that $X=K\cup L$.
If every nondegenerate subcontinuum of $X$ is decomposable,
then $X$ is {\it hereditarily decomposable}.
The continuum $X$ is {\it indecomposable} if it is not decomposable.
If each subcontinuum of $X$ is indecomposable, then $X$ is
{\it hereditarily decomposable}.
If $X$ is a continuum, a subcontinuum $K$ of $X$ is {\it terminal}
provided that for each subcontinuum $L$ of $X$ satisfying that
$K\cap L\not=\emptyset$, we have that either $K\subset L$ or
$L\subset K$.

A continuum $X$ is
{\it chainable} provided that for each $\varepsilon>0$, there exists a finite
open cover $\{U_1,\ldots, U_n\}$ such that $U_j\cap U_k\not=\emptyset$ if and only
if $|j-k|\leq 1$, and $\diam(U_j)<\varepsilon$, for all $j\in\{1,\ldots,n\}$. 
The continuum $X$ is {\it circularly chainable} if for each $\varepsilon>0$, there exists a finite
open cover $\{U_1,\ldots, U_n\}$ such that $U_j\cap U_k\not=\emptyset$ if and only
if $|j-k|\leq 1$ or $\{j,k\}=\{1,n\}$, and $\diam(U_j)<\varepsilon$, for every $j\in\{1,\ldots,n\}$. 

The {\it pseudo-arc} is a hereditarily indecomposable chainable continuum $P$. 
R H Bing shows that the pseudo-arc is homogeneous \cite[Theorem 13]{Bing}.
He also proves that any two pseudo-arcs are homeomorphic \cite[Theorem 1]{Bi}. 
The {\it circle of pseudo-arcs} is a planar homogeneous continuum $Z$ 
admitting a continuous decomposition $\G$ into pseudo-arcs such that $Z/\G$ is a simple closed curve.
R H Bing and F. B. Jones prove that any two circle of pseudo-arcs are homeomorphic \cite[p. 179]{BJ}.
The {\it arc of pseudo-arcs} is a chainable continuum with a continuous decomposition into pseudo-arcs such that the
decomposition space is an arc. R H Bing and F. B. Jones show that the arc of pseudo-arcs is unique and
proved some theorems on extending homeomorphisms on arcs of pseudo-arcs \cite{BJ}.
Note that, by construction, the arc of pseudo-arcs is a saturated subcontinuum of the circle of pseudo-arcs. 
A construction of the circle of pseudo-arcs and the arc of pseudo-arcs and, in general,
continuous curves of pseudo-arcs may be found in \cite{L}.

A continuum $X$ is {\it irreducible} if there exist two points
$p$ and $q$ of $X$ such that no proper subcontinuum of $X$ contains
both $p$ and $q$. A continuum $X$ is of {\it type $\lambda$}
provided that $X$ is irreducible and each indecomposable subcontinuum
of $X$ has empty interior. By \cite[Theorem 10, p. 15]{T}, a continuum $X$
is of type $\lambda$ if and only if it admits the finest monotone upper semicontinuous
decomposition $\G$ such that $X/\G$ is an arc. Each element of $\G$
is called a {\it layer} of $X$. Following \cite{MO}, we say that
a continuum $X$ of type $\lambda$ for which $\G$ is continuous, 
is a {\it continuously irreducible continuum}.

A continuum $X$ is a {$\theta$-continuum} 
({$\theta_n$-continuum}, for some positive integer $n$) if for each
subcontinuum $K$ of $X$, we have that $X\setminus K$ only has
a finite number of components ($X\setminus K$ has at most
$n$ components). Following \cite{Vo2}, we say that
$\theta$-continuum ($\theta_n$-continuum) is of {\it type $A$}
provided that it admits a monotone upper semicontinuous
decomposition $\D$ whose quotient space is a finite graph, and it
is of {\it type $A'$} if, in addition, the elements of the decomposition
have empty interior. A $\theta$-continuum ($\theta_n$-continuum) 
of type $A'$ for which the decomposition $\D$ is continuous is a
{\it continuously type $A'$ $\theta$-continuum $($$\theta_n$-continuum$)$}.

A {\it generalized continuum} is a locally compact, connected, metric space.
A generalized continuum $X$ is an {\it exhausted $\sigma$-continuum} 
provided that there
exists a sequence $\{K_n\}_{n=1}^\infty$ of subcontinua of $X$ 
such that $K_n\subset Int(K_{n+1})$,
for all $n\in\NN$, and $X=\bigcup_{n=1}^\infty K_n$. The sequence 
$\{K_n\}_{n=1}^\infty$ is
an {\it exhaustive sequence} of subcontinua of $X$.

Let $X$ be a continuum. We consider the following
{\it hyperspaces} of $X$:
$$2^X=\{A\subset X\ |\ A\ \hbox{is closed and nonempty}\};$$
$$\C_1(X)=\{A\in 2^X\ |\ A\ \hbox{is connected}\};$$
$$\FF_1(X)=\{\{x\}\ |\ x\in X\}.$$
$2^X$ is topologised with the Hausdorff metric, $\HH$, given by:
$$\HH(A,B)=\inf\{\varepsilon>0\ |\ A\subset\V_\varepsilon(B)\ \hbox{and}\ B\subset\V_\varepsilon(A)\}.$$
Given a finite collection, $U_1,\ldots,U_m$, of nonempty subsets of $X$, we
define:
$$\langle U_1,\ldots,U_m\rangle\!=\!\{A\in 2^X\ \! |\ \! A\subset\bigcup_{k=1}^mU_k\ \hbox{\rm and
$A\cap U_k\not=\emptyset$, for all $k\in\{1,\ldots,m\}$}\}.$$
It is known that the family of all subsets of $2^X$ of the form $\langle U_1,\ldots,U_m\rangle$, where
$U_1,\ldots, U_m$ are nonempty open subsets of $X$,
form a basis for a topology for $2^X$ called {\it Vietoris Topology} \cite[(0.11)]{N1}, and
that the Vietoris Topology and the Topology induced by the Hausdorff metric coincide 
\cite[(0.13)]{N1}. If $A$ and $B$ are two elements of $2^X$, an {\it order arc} from $A$ to $B$
is a one-to-one map $\alpha\colon [0,1]\to 2^X$ such that $\alpha(0)=A$, $\alpha(1)=B$ and
$\alpha(s)\subsetneq\alpha(t)$ if $s<t$. 

A continuum $X$ has the {\it property of Kelley at a point $p\in X$} 
provided that for each $\varepsilon>0$, there exists a $\delta>0$ such
that if $A$ is a subcontinuum of $X$ and $p\in A$, then for each point
$q\in X$ such that $d(p,q)<\delta$, there exists a subcontinuum $B_q$ 
of $X$ such that $q\in B_q$ and $\HH(A,B_q)<\varepsilon$.
The number $\delta$ in the definition is called a {\it Kelley number} for
$(X,\varepsilon,p)$. A continuum $X$ has the {\it property of Kelley} if
$X$ has the property of Kelley at each point.

A continuum $X$ has the {\it property of Kelley weakly} provided that 
there exists a dense subset $\A$ of $\C_1(X)$ such that $X$ has
the property of Kelley at some point of each $A\in\A$. This concept
is introduced in \cite{MN}. An example of a continuum with the
property of Kelley weakly that does not have the property of Kelley
is presented in \cite[Example 2.1]{MN}.

Given a continuum $X$, we define Professor Jones' 
{\it Set Function $\T$} as follows: if $A$ is a subset of $X$, then
\begin{multline*} 
\T(A)=X\setminus\{x\in X\ |\ \hbox{\rm there exists a subcontinuum}\\
\hbox{\rm $W$ of $X$ such that $x\in Int(W)\subset W\subset$ $X\setminus A$}\}.
\end{multline*}
Let us observe that for each 
subset $A$ of $X$, $\T(A)$ is a closed subset of $X$ and $A\subset\T(A)$. 
When there is a possibility of
confusion, we add a subscript to $\T$ to indicate the
continuum we are considering.
Note that, if $W$ is a subcontinuum of $X$, then $\T(W)$ is a subcontinuum
of $X$ \cite[Theorem 2.1.27]{MS}. Also, by \cite[Corollary 2.1.14]{MS},
$\T(\emptyset)=\emptyset$. A nonempty closed subset $A$ of $X$ is a {\it $\T$-closed
set} if $\T(A)=A$. The family of $\T$-closed subsets of $X$ is denoted by $\mathfrak{T}(X)$,
and $\mathfrak{T}_\C(X)=\mathfrak{T}(X)\cap\C_1(X)$.
The set function $\T$ is {\it idempotent} if 
$\T^2(A)=\T(A)$ for all subsets $A$ of $X$. Now, $\T$ is {\it idempotent on closed sets}
({\it idempotent on continua}, {\it idempotent on singletons}, respectively) if for each nonempty
closed subset (continuum, singleton, respectively) $A$ of $X$, $\T^2(A)=\T(A)$.
Since the restriction of $\T$ to $2^X$ sends $2^X$ into $2^X$, and $2^X$ we
may ask if $\T\colon 2^X\to 2^X$ is continuous. Hence, we say that
{\it $\T$ is continuous} if it restriction to $2^X$ is continuous. Also,
since $\T(\C_1(X))\subset\C_1(X)$, we say that {\it $\T$ is continuous
on continua} provided that $\T\colon\C_1(X)\to\C_1(X)$ is continuous.

Let $X$ be a continuum and let $A$ be a subcontinuum of $X$. Then
\smallskip

$\bullet$ $A$ is a {\it subcontinuum of colocal connectedness} 
of $X$ if for each
open subset $U$ of $X$ containing $A$, there exists an open subset
$V$ of $X$ such that $A\subset V\subset U$ and $X\setminus V$ is
connected.
\smallskip

$\bullet$ $X\setminus A$ is {\it strongly continuumwise connected} provided that
for each pair of points $x_1$ and $x_2$ of $X\setminus A$, there exists
a subcontinuum $W$ of $X\setminus A$ such that $\{x_1,x_2\}\subset Int(W)$.
\smallskip

$\bullet$ $X\setminus A$ is {\it continuumwise connected} if 
for each pair of points $x_1$ and $x_2$ of $X\setminus A$, there exists
a subcontinuum $W$ of $X\setminus A$ such that $\{x_1,x_2\}\subset W$.
\smallskip

$\bullet$ $A$ is a {\it nonblock} subcontinuum of $X$ provided that 
there exists an increasing
sequence $\{L_n\}_{n=1}^\infty$ of subcontinua of 
$X\setminus A$ such that 
$\bigcup_{n=1}^\infty L_n$ is dense in $X\setminus A$.
\smallskip

$\bullet$ $A$ is a {\it shore} subcontinuum of $X$, 
if for every
$\varepsilon>0$, there exists a subcontinuum $W$ of 
$X\setminus A$ such that 
each $\varepsilon$-ball in $X$ intersects $W$.
\smallskip

$\bullet$ $A$ is {\it not a strong centre} in $X$ provided that for every pair of open
subsets $U$ and $V$ of $X$, then there exists a subcontinuum $W$ of
$X$ such that $W\cap U\not=\emptyset$, $W\cap V\not=\emptyset$ and
$W\cap A=\emptyset$.

\begin{remark}
{\rm Note that when $A$ is a nonblock, a shore subcontinuum, and
$A$ is not a strong centre, we require the 
interior of $A$ to be empty.}
\end{remark}

Next lemma follows from definitions.

\begin{lemma}\label{lemyhfgt5}
Let $X$ be a continuum and let $A$ be a subcontinuum of $X$
such that $Int(A)=\emptyset$. Consider the following statements:
\smallskip

$(1)$  $A$ is a continuum of colocal connectedness of $X$;
\smallskip

 $(2)$ $X\setminus A$ is strongly continuumwise connected;
 \smallskip
 
$(3)$ $X\setminus A$ is continuumwise connected;
\smallskip

$(4)$ $A$ is a nonblock subcontinuum of $X$;
\smallskip

$(5)$ $A$ is a shore subcontinuum of $X$;
\smallskip

$(6)$ $A$ is not a strong centre in $X$.
\smallskip

\noindent Then item $(j)$ implies item $(k)$, whenever $j\in\{1,2,3,4,5\}$,
$k\in\{2,3,4,5,6\}$, and $j<k$.
\end{lemma}

\begin{theorem}\label{lemma8uhgvfr5}
Let $X$ be a continuum and let $A$ be a subcontinuum of $X$. Then 
$A$ is a subcontinuum of colocal connectedness if and only if $X\setminus A$ 
is strongly continuumwise connected.
\end{theorem}

{\bf Proof.} By Lemma~\ref{lemyhfgt5}, if $A$ is a subcontinuum of colocal 
connectedness of $X$, 
then $X\setminus A$ is strongly continuumwise connected. Suppose that 
$X\setminus A$ is strongly continuumwise connected, we prove that $A$ is 
a subcontinuum of colocal connectedness of $X$. Let $V$ be an open subset of 
$X$ such that $A\subseteq V$. By \cite[Theorem 1.4.75]{MS}, there exists a 
sequence of subcontinua $\{W_n\}_{n=1}^\infty$ of $X$ such that 
$X\setminus A=\bigcup_{n=1}^\infty W_n$ and $W_n\subseteq Int(W_{n+1})$,  
for each $n\in\NN$. Note that 
$A=\bigcap_{n=1}^\infty (X\setminus W_n)=\bigcap_{n=1}^\infty (X\setminus Int(W_n))$,  
where $X\setminus Int(W_n)$ is compact, and $X\setminus Int(W_{n+1})\subseteq X\setminus Int(W_n)$ 
for all $n\in\NN$. Thus, there exists $k\in\NN$ such that 
$X\setminus Int(W_k)\subseteq V$ \cite[Lemma 1.6.7]{MT}. Let $U=X\setminus W_k$. Observe that 
$A\subseteq U\subseteq V$ and $X\setminus U$ is connected. Therefore, $A$ is a subcontinuum of colocal connectedness.\par\Q

Note that, by Theorem~\ref{lemma8uhgvfr5}, we know
that $(1)$ and $(2)$ of Lemma~\ref{lemyhfgt5} are equivalent. Also,
in \cite{Bobok}, the reader can find examples of continua where the 
authors show that the other reverse implications in Lemma~\ref{lemyhfgt5} 
do not hold for degenerate subcontinua.

\begin{theorem}\label{strongcontconniffcontconn}
Let $X$ be a continuum with the property of Kelley and let $A$ be a 
subcontinuum of $X$. Then the following are equivalent:
\smallskip

$(1)$ $A$ is a subcontinuum of colocal connectedness of $X$.
\smallskip

$(2)$ $X\setminus A$ is strongly continuumwise connected.
\smallskip

$(3)$ $X\setminus A$ is continuumwise connected.
\end{theorem}

{\bf Proof.} By Theorem~\ref{lemma8uhgvfr5}, we have that $(1)$
and $(2)$ are equivalent. 
By Lemma~\ref{lemyhfgt5}, $(2)$ implies $(3)$. Assume that
$X\setminus A$ is continuumwise connected. Let $x_1$ and $x_2$ be two
points of $X\setminus A$. Since $X\setminus A$ is continuumwise connected,
there exists a subcontinuum $K$ of $X\setminus A$ such that 
$\{x_1,x_2\}\subset K$.
Let $n\in\NN$ and let $W_n$ be the component of $X\setminus\V_{1\over n}(A)$
containing $K$. Note that $W_n\subset W_{n+1}$, for every $n\in\NN$.
Let $x_3\in X\setminus A$. Then there exists a subcontinuum
$L$ of $X\setminus A$ such that $\{x_2,x_3\}\subset L$. Hence, $K\cup L$ is
a subcontinuum of $X\setminus A$. Thus, there exists $m\in\NN$ such that
$K\cup L\subset W_m$. Hence, $X\setminus A=\bigcup_{n=1}^\infty W_n$.
Since $X\setminus A$ is open, by the Baire Category Theorem 
\cite[Theorem 25.3]{Will},
there exists $k\in\NN$ such that $Int(W_k)\not=\emptyset$. Let $r>\max\{m,k\}$.
Then $K\cup L\subset W_r$ and $Int(W_r)\not=\emptyset$. Since
$W_r\cap A=\emptyset$, by \cite[Corollary 1.6.21]{MS}, there 
exists a subcontinuum
$M$ of $X$ such that $W_r\subset Int(M)$ and $M\cap A=\emptyset$.
Therefore, $X\setminus A$ is strongly continuumwise connected,
and $(3)$ implies $(2)$.\par\Q

The following is \cite[Lemma (Bing)]{Vo2}.

\begin{theorem}\label{LemmaBing}
If $X$ is a $\theta_1$-continuum,
then each proper subcontinuum of $X$ is not a strong centre in $X$.
\end{theorem}

\begin{theorem}\label{nonblockset1}
Let $X$ and $Y$ continua and let $f\colon X\tto Y$ be a monotone open map. If $B$ is 
a nonblock subcontinuum of $Y$, then $f^{-1}(B)$ is a nonblock subcontinuum of $X$.
\end{theorem}

{\bf Proof.} Since $B$ is a nonblock subcontinuum of $Y$, there exists an increasing sequence $\{K_n\}_{n=1}^\infty$ of subcontinua
of $Y\setminus B$ such that $\bigcup_{n=1}^\infty K_n$ is dense in $Y\setminus B$. 
Since $f$ is monotone, $\{f^{-1}(K_n)\}_{n=1}^\infty$ is an increasing sequence of subcontinua of $X$. 
Note that $f^{-1}(B)$ is a continuum and $f^{-1}(B)\cap f^{-1}(K_n)=\emptyset$ for each $n\in\NN$. 

We show that $\bigcup_{n=1}^\infty f^{-1}(K_n)$ is dense in $X$. Let $U$ be a nonempty open subset of $X$. 
Since $f$ is open, $f(U)$ is a nonempty open subset of $Y$. Hence, there exists $m\in\mathbb N$ such that $f(U)\cap K_m\neq\emptyset$. 
Thus, $U\cap f^{-1}(K_m)\neq\emptyset$ and $\bigcup_{n=1}^\infty f^{-1}(K_n)$ is dense in $X$. 
Therefore, $f^{-1}(B)$ is a nonblock subcontinuum of $X$.\par\Q

The following lemma is useful in the proof of Theorem~\ref{ATclosedequiv}.

\begin{lemma}\label{lem7yg6}
Let $X$ be a continuum and let $A$ be a nonempty closed subset of $X$. 
Then the following are equivalent:
\smallskip

$(1)$ For each open subset $U$ of $X$ such that $A\subseteq U$ 
there exists an open subset
$V$ of $X$ such that $A\subseteq V\subseteq U$ and $X\setminus V$ has 
finitely many components;
\smallskip
        
$(2)$ For each open subset $U$ of $X$ such that $A\subseteq U$ 
there exists an open subset $V$ of $X$ such that 
$A\subseteq V\subseteq U$ and $X\setminus V =N_1\cup \cdots \cup N_k$, where 
$N_1,\ldots,N_k$ are pairwise disjoint subcontinua of $X$, 
and $Int(N_i)\neq\emptyset$ for each $i\in\{1,\ldots,k\}$.
\end{lemma}

\textbf{Proof.} It is clear that $(2)$ implies $(1)$. Suppose $(1)$ 
and let $U$ be an open subset of 
$X$ such that $A\subseteq U$. Let $U'$ be an open set such that 
$A\subseteq U'\subseteq Cl(U')\subseteq U$. 
By $(1)$, there exists an open subset $W$ of $X$ such that 
$A\subseteq W\subseteq U'$ and $X\setminus W$ has 
only finitely many components. Let $M_1,\ldots,M_l$ be the components 
of $X\setminus W$. Observe 
that $Int(M_i)=\emptyset$ if and only if $M_i\subseteq Cl(W)\setminus W$. Let 
$V=W\cup\{M_i\ |\ Int(M_i)=\emptyset\}$. Hence, we have that $V$ is an 
open subset of $X$, 
$V\subseteq Cl(U')\subseteq U$ and $X\setminus V=N_1\cup \cdots \cup N_k$, 
where $N_1,\ldots,N_k$ are pairwise disjoint subcontinua of $X$, 
and $Int(N_i)\neq\emptyset$, 
for each $i\in\{1,\ldots,k\}$.\par\Q

\begin{lemma}\label{lemmat}
Let $X$ be an irreducible continuum and let $p$ be a point of 
irreducibility Then 
\smallskip

$(1)$ If $Cl(X\setminus \kappa(p))$ has empty interior, then 
$Cl(X\setminus \kappa(p))=X\setminus \kappa(p)$.
\smallskip

$(2)$ If $Cl(X\setminus \kappa(p))$ has nonempty interior, then 
$Cl(X\setminus \kappa(p))$ is indecomposable.
\end{lemma}

\textbf{Proof.} To prove $(1)$, suppose that 
$z\in Cl(X\setminus\kappa(p))\setminus (X\setminus\kappa(p))$. 
Then there exists a proper subcontinuum $L$ of $X$ 
such that $\{z,p\}\subseteq L$. Observe that $X=Cl(X\setminus\kappa(p))\cup L$. 
Hence, If $U=X\setminus L$, then $U$ is a nonempty open subset of $X$ such that 
$U\subseteq Cl(X\setminus\kappa(p))$, a contradiction. Therefore, 
$Cl(X\setminus \kappa(p))=X\setminus \kappa(p)$.

\bigskip

To show $(2)$, assume that $Cl(X\setminus\kappa(p))$ has nonempty interior and
let $K=Cl(X\setminus\kappa(p))$. We consider two cases:
\medskip

\textbf{Case (1).} $X=K$. 
\smallskip

Suppose that there exist two proper subcontinua $A$ and $B$ of $X$ 
such that $X=A\cup B$. Note that $A\cap (X\setminus \kappa(p))\neq\emptyset$ 
and $B\cap (X\setminus \kappa(p))\neq\emptyset$. Suppose that $p\in A$. Hence, 
$A=X$, a contradiction. Similarly, we obtain a contradiction if $p\in B$. 
Therefore, $X=K$ is indecomposable.
\medskip

\textbf{Case (2).} $X\setminus K\neq\emptyset$. 
\smallskip

Let $L=Cl(X\setminus K)$. We know that $K$ is a proper subcontinuum of 
$X$ such that $p\in L$. Suppose that $K$ is decomposable; i.e., there exist 
two proper subcontinua $A$ and $B$ of $K$ such that $K=A\cup B$. 
Since $L\cap K\neq\emptyset$, 
either $L\cap A\neq\emptyset$ or $L\cap B\neq\emptyset$. Suppose that 
$L\cap A\neq\emptyset$. Since $X\setminus\kappa(p)\neq\emptyset$, $X=L\cup A$ 
and $K=A$, a contradiction.\par\Q

\begin{theorem}\label{theo0}
Let $X$ and $Y$ be chainable continua. Let $A$ and $B$ be subcontinua 
of $X\times Y$. 
If $\pi_X(A)=X$ and $\pi_Y(B)=Y$, then $A\cap B\neq\emptyset$.
\end{theorem}

{\bf Proof.} Let $A$ and $B$ be subcontinua of $X\times Y$ 
such that  $\pi_X(A)=X$ 
and  $\pi_Y(B)=Y$. Suppose that  $A\cap B=\emptyset$.

Let $d_X$ and $d_Y$ be metrics on $X$ and $Y$, respectively. 
We consider a metric 
$\rho$ on $X\times Y$ defined by 
$$\rho((x_1,y_1),(x_2,y_2))=\max\{d_X(x_1,x_2),d_Y(y_1,y_2)\}.$$ 
Since $A$ and $B$ are disjoint, there exists $\varepsilon>0$ such that 
\begin{equation}\label{ec00}
  \min\{\rho(a,b)\ |\ a\in A \text{ and }b\in B\}>\varepsilon.  
\end{equation}

Since $X$ and $Y$ are chainable, there exist surjective $\varepsilon$-maps
$f\colon X\tto [0,1]$ 
and $g\colon Y\tto [0,1]$. Let $G\colon X\times Y\tto [0,1]^2$ 
be defined for each $(x,y)\in X\times Y$ by 
$$G((x,y))=(f(x),g(y)).$$

Note that, since $\pi_X(A)=X$, $\pi_1(G(A))=[0,1]$. Hence, 
$G(A)\cap \{0\}\times [0,1]\neq\emptyset$ and 
$G(A)\cap \{1\}\times [0,1]\neq\emptyset$. 
Similarly, since $\pi_Y(B)=Y$, we have that 
$G(B)\cap ([0,1]\times \{0\})\neq\emptyset$ and 
$G(B)\cap ([0,1]\times \{1\})\neq\emptyset$. Therefore, by 
\cite[Theorem 130, p. 158]{Moore}, 
$G(A)\cap G(B)\neq\emptyset$. Let $(a_1,a_2)\in A$ and $(b_1,b_2)\in B$ 
be such that 
$G((a_1,a_2))=G((b_1,b_2))$. Thus, $f(a_1)=f(b_1)$ and $g(a_2)=g(b_2)$. 
Since $f$ and $g$ 
are $\varepsilon$-maps, $d_X(a_1,b_1)<\varepsilon$ and 
$d_Y(a_2,b_2)<\varepsilon$. 
Thus, $\rho ((a_1,a_2),(b_1,b_2))<\varepsilon$. We contradict (\ref{ec00}). 
Therefore, $A\cap B\neq\emptyset$. \par\Q

Since each proper subcontinuum of a circularly chainable 
continuum is chainable \cite[Corollary 2.1.45]{MT},
as a consequence of Theorem~\ref{theo0}, we have:

\begin{corollary}
Let $X$ and $Y$ be circularly chainable continua, and let 
$A$ and $B$ proper subcontinua
of $X$ and $Y$, respectively. If $C$ and $D$ are subcontinua 
of $A\times B$ such that
$\pi_X(C)=A$ and $\pi_Y(D)=B$, then $C\cap D\neq\emptyset$.
\end{corollary}

\section{Properties of $\Delta_X$}\label{Diagonal}

For each continuum $X$, we denote the diagonal of $X^2$ by 
$\Delta_X=\{(x,x)\ |\ x\in X\}$. 

In this section we present results of independent interest
regarding properties of the diagonal of a continuum
which are used in section~\ref{DeltaTclosed}. 

\begin{theorem}\label{theodescdiagnequ}
If $X$ is a decomposable continuum, then there exists a subcontinuum
$L$ of $X^2\setminus\Delta_X$ such that $Int_{X^2}(L)\not=\emptyset$. 
\end{theorem}
\textbf{Proof.} Let $A$ and $B$ be proper subcontinua such that $X=A\cup B$. 
Let $U$ and $V$ be open subsets of $X$ such that $Cl(U)\cap B=\emptyset$ and 
$Cl(V)\cap A=\emptyset$. Let $L=(Cl(U)\times B)\cup (A\times Cl(V))$. 
Then we have 
that $L$ is a subcontinuum of $X^2\setminus\Delta_X$ such that 
$U\times V\subseteq L$.\par\Q

\begin{theorem}\label{theta1cont}
If $X$ is a $\theta_1$-continuum, then $\Delta_X$ is not 
a strong centre in $X^2$.
\end{theorem}

{\bf Proof.} Let $\U$ and $\V$ open subsets of $X^2$. We may find
pairwise disjoint open subsets $U_1$, $U_2$, $V_1$ and $V_2$ of $X$
such that $U_1\times U_2\subset\U$ and $V_1\times V_2\subset\V$.
Let $L_1$ be a subcontinuum of $X$ such that 
$L_1\cap U_1\not=\emptyset$ and $L_1\cap V_1\not=\emptyset$. 
By Theorem~\ref{LemmaBing}, there exists a subcontinuum
$L_2$ of $X\setminus L_1$ such that $L_2\cap U_2\not=\emptyset$ 
and $L_2\cap V_2\not=\emptyset$. Then $W=L_1\times L_2$
is a subcontinuum of $X^2\setminus\Delta_X$ such that
$W\cap (U_1\times U_2)\not=\emptyset$ and 
$W\cap (V_1\times V_2)\not=\emptyset$. Therefore,
$\Delta_X$ is not a strong centre in $X^2$.\par\Q

\begin{theorem}\label{typelambda}
Let $X$ be a chainable continuum. If $X$ is of the type $\lambda$, 
then $\Delta_X$ is a strong centre.
\end{theorem}

\textbf{Proof.} Let $m\colon X\tto [0,1]$ be the monotone map such that 
$\{m^{-1}(t)\ |\ t\in [0,1]\}$ is the finest admissible decomposition of 
$X$. Let $U=m^{-1}([0,\frac{1}{4}))\times m^{-1}((\frac{3}{4},1])$ and let 
$V=m^{-1}((\frac{3}{4},1])\times m^{-1}([0,\frac{1}{4}))$. Observe that
$U$ and $V$ are open subsets of $X^2$.

We show that if $L$ is a subcontinuum of $X^2$ such that 
$L\cap U\neq\emptyset$ and $L\cap V\neq\emptyset$, then 
$L\cap \Delta_X\neq\emptyset$. Suppose that $L$ is a subcontinuum 
of $X$ such that $L\cap U\neq\emptyset$, $L\cap V\neq\emptyset$ and 
$L\cap \Delta_X=\emptyset$. Let 
$A=m^{-1}([0,\frac{1}{4}])\times m^{-1}([\frac{3}{4},1])$ and let 
$B=m^{-1}([\frac{3}{4},1])\times m^{-1}([0,\frac{1}{4}])$. Note that $A$ 
and $B$ are proper subcontinua of $X^2\setminus\Delta_X$. Since 
$U\subseteq A$ and $V\subseteq B$, $L\cap A\neq\emptyset$ and 
$L\cap B\neq\emptyset$. Hence, $L\cup A\cup B$ is a subcontinuum of
$X^2\setminus\Delta_X$. Observe that $\pi_1(L\cup A\cup B)=
\pi_2(L\cup A\cup B)=X$ 
and $\Delta_X\cap (L\cup A\cup B)=\emptyset$. Since $X$ is chainable, we obtain
a contradiction to Theorem \ref{theo0}. Therefore, $\Delta_X$ is a 
strong centre in $X^2$.\par\Q

Since every Elsa continuum is a type $\lambda$ continuum that is chainable, 
the following corollary generalizes \cite[Proposition 26]{IMMM2}.

\begin{corollary}
Let $X$ be chainable continuum of type $\lambda$ different from an arc. Then
$X^2\setminus\Delta_X$ is 
connected, but $\Delta_X$ is a strong centre.
\end{corollary}

{\bf Proof.} By \cite[Theorem~3]{Kat}, $X^2\setminus\Delta_X$ is connected 
and, by Theorem \ref{typelambda}, $\Delta_X$ is a strong centre. \par\Q

\begin{theorem}\label{chaintermsubcont}
Let $X$ be a continuum and let $K$ be a terminal subcontinuum of $X$. 
If $K$ is chainable, then $X^2\setminus \Delta_X$ is not continuumwise connected.
\end{theorem}

\textbf{Proof.} Since $K$ is a chainable continuum, by \cite[Theorem 12.5]{N2},
there exist two elements $a$ and $b$ in $K$ be such that $K$ is irreducible 
between $a$ and $b$. We show that there does
not exist a subcontinuum $L$ of $X^2$ such that $\{(a,b),(b,a)\}\subseteq L$ 
and $L\cap \Delta_X=\emptyset$. Suppose there exists a subcontinuum $L$ 
of $X^2$ such that $\{(a,b),(b,a)\}\subseteq L$. Since 
$\{a,b\}\subseteq \pi_1(L)\cap \pi_2(L)$, we 
have that $K\subseteq \pi_1(L)\cap\pi_2(L)$. Let $\alpha\colon [0,1] \to\C_1(X^2)$ 
be an order arc such that $\alpha(0)=\{(a,b)\}$ and $\alpha(1)=L$
\cite[Theorem (1.8)]{N1}. Let 
$$t_1=\max\{t\in [0,1]\ |\ \pi_1(\alpha(t))\subseteq K\}\text{ and }t_2=
\max\{t\in [0,1]\ |\ \pi_2(\alpha(t))\subseteq K\}.$$ 
Note that, since $K$ is a terminal subcontinuum of $X$, 
$\pi_1(\alpha(t_1))=K$ and $\pi_2(\alpha(t_2))=K$. 
Suppose that $t_1\leq t_2$. Let $R=\alpha(t_1)$. Observe that $R$ and $\Delta_K$ 
are subcontinua of $K^2$ such that $\pi_1(R)=K$ and $\pi_2(\Delta_K)=K$. Thus, 
by Theorem \ref{theo0}, $R\cap \Delta_K\neq\emptyset$. Since $R\subseteq L$ and 
$\Delta_K\subseteq \Delta_X$, $L\cap\Delta_X\neq\emptyset$. \par\Q

Since every subcontinuum of a chainable continuum is chainable
\cite[Theorem (9.C.4)]{Chris}, next corollary follows from Theorem~\ref{chaintermsubcont}.

\begin{corollary}
Let $X$ be a chainable continuum. If there exists a nondegenerate terminal 
subcontinuum $K$ of $X$, then $X^2\setminus \Delta_X$ is not continuumwise connected.    
\end{corollary}


If $X$ is a one-dimensional continuum, then Wayne Lewis \cite {L} 
constructed what is known as the
{\it continuum of pseudo-arcs of $X$}, which is a one-dimensional 
continuum admitting a continuous
terminal decomposition, $\G$, into pseudo-arcs whose quotient 
space is $X$. Hence, since the
pseudo-arc is a chainable continuum \cite[Definition 7]{Moise}, as a
consequence of Theorem~\ref{chaintermsubcont}, we have:

\begin{corollary}
Let $X$ be a one-dimensional continuum. If $\widehat X$ is the 
continuum of pseudo-arcs of $X$, then
$\widehat X^2\setminus \Delta_{\widehat X}$ is not continuumwise connected.
\end{corollary}

\begin{theorem}\label{nonblockset}
Let $X$ and $Y$ be continua, where $Y$ is locally connected, and let
$f\colon X\tto Y$ be a monotone open map. If $Y$ is not an arc, then
$\Delta_X$ is a nonblock subcontinuum of $X^2$.
\end{theorem}

{\bf Proof.} Since $Y$ is a locally connected continuum that is not an arc, 
by \cite[Corollary 6]{IMMM2}, $\Delta_Y$ is a nonblock subset of $Y^2$. Let 
$f\times f\colon X^2\to Y^2$ be the map defined by $(f\times f)((x,x'))=(f(x),f(x'))$ 
for each $(x,x')\in X^2$. Note that $f\times f$ is a monotone open map. By 
Theorem~\ref{nonblockset1}, $(f\times f)^{-1}(\Delta_Y)$ is a nonblock subcontinuum 
of $X^2$. Since $\Delta_X\subseteq (f\times f)^{-1}(\Delta_Y)$, $\Delta_X$ is a 
nonblock subcontinuum of $X^2$.
\par\Q


%

\begin{remark}\label{f_2}
Let $f_2\colon X^2\tto\FF_2(X)$ be given by $f_2((x_1,x_2))=\{x_1,x_2\}$. Then $f_2$ is an open map \cite[Lemma 9]{M0}.
Note that $f_2(\Delta_X)=\FF_1(X)$ and $f_2^{-1}(\FF_1(X))=\Delta_X$. If $\M$ is a subcontinuum of $\FF_2(X)\setminus\FF_1(X)$,
then $f_2^{-1}(\M)=M\cup M^*$, where $M^*=\{(x_2,x_1)\in X^2\ |\ (x_1,x_2)\in M\}$ and, in this case, 
$M$ and $M^*$ are the components of $f_2^{-1}(\M)$  \cite[Lemma~12]{M0}.
\end{remark}

\begin{theorem}\label{scimpliesnonblock}
Let $X$ be a continuum. If $\Delta_X$ is not a strong centre in $X^2$, 
then $\Delta_X$ is a nonblock subcontinuum of $X^2$. 
\end{theorem}

\textbf{Proof.} If $X$ is indecomposable, then $\Delta_X$ is a 
nonblock subcontinuum 
of $X^2$ \cite[Theorem 28]{IMMM2}. Suppose now that $X$ is decomposable. 
By Theorem \ref{theodescdiagnequ}, there exists a subcontinuum $L$ of 
$X^2\setminus \Delta_X$ such that $Int(L)\neq\emptyset$. Let $U$ and $V$ 
be disjoint open subsets of $X$ such that $U\times V\subseteq L$. Since 
$\Delta_X$ is not a strong centre, there exists a continuum $J$ such that 
$J\cap (U\times V)\neq\emptyset$, $J\cap (V\times U)\neq\emptyset$ and 
$J\cap \Delta_X=\emptyset$.
Let $f_2\colon X^2\tto\FF_2(X)$ be the map given in Remark~\ref{f_2}.
Note that by \cite[Lemma 9]{M0}, $f_2$ is open.
Thus, $f_2(U\times V)$ is an open subset of $X^2\setminus\Delta_X$. By 
\cite[Theorem 3.5]{MM}, there exists a sequence of subcontinua 
$\{\M_n\}_{n=1}^\infty$ of $\FF_2(X)\setminus\FF_1(X)$ such that 
$\M_n\subseteq\M_{n+1}$ and $\bigcup_{n=1}^\infty\M_n$ is dense in 
$\FF_2(X)\setminus\FF_1(X)$. Without loss of generality, we may assume that 
$\M_1\cap f_2(U\times V)\neq\emptyset$. 
Hence, by \cite[Lemma 12]{M0}, $f_2^{-1}(\M_n)=M_n\cup M^*_n$, 
where $M_n$ and $M^*_n$ 
are subcontinua of $X^2\setminus\Delta_X$, $M_n\cap (U\times V)\neq\emptyset$ 
and $M^*_n\cap (V\times U)\neq\emptyset$. Let $L^*=\{(x,y)\ |\ (y,x)\in L\}$. 
Note that $V\times U\subseteq L^*$. For each $n\in\NN$, let 
$$K_n=M_n\cup L\cup J\cup L^*\cup M^*_n.$$ 
Observe that $K_n$ is a subcontinuum of $X^2$ such that 
$f_2^{-1}(\M_n)\subseteq K_n$ 
and $K_n\cap\Delta_X=\emptyset$. Since $\bigcup_{n=1}^\infty\M_n$ 
is dense in $\FF_2(X)$, we obtain that $\bigcup_{n=1}^\infty K_n$ is dense in 
$X^2\setminus\Delta_X$. Therefore, $\Delta_X$ is a nonblock 
subcontinuum of $X^2$. \par\Q

The next corollary gives a negative answer to 
\cite[Questions 35 and 36]{IMMM2}.
It also gives a negative answer to \cite[Question 38]{IMMM2} for 
$i\in\{4,5\}$.

\begin{corollary}\label{scimpliesnonblock2}
Let $X$ be a continuum. Then the following are equivalent:
\smallskip

$(1)$ $\Delta_X$ is a nonblock subcontinuum of $X^2$.
\smallskip

$(2)$ $\Delta_X$ is a shore subcontinuum of $X^2$.
\smallskip

$(3)$ $\Delta_X$ is not a strong centre in $X^2$.
\end{corollary}

{\bf Proof.} By Lemma~\ref{lemyhfgt5}, we have that $(1)$ implies $(2)$ and
$(2)$ implies $(3)$. Now, by Theorem~\ref{scimpliesnonblock}, we obtain that 
$(3)$ implies $(1)$.\par\Q

The next proposition gives a negative answer to \cite[Question 1.2]{MM}. 

\begin{proposition}\label{examnegMM}
    There exists a nonlocally connected continuum $X$ such that $\FF_1(X)$ is a continuum of colocal connectedness of $\FF_2(X)$. 
\end{proposition}

{\bf Proof.} We denote $xy=\{tx+(1-t)y\ |\ t\in [0,1]\}$, for any $x,y\in\RR^2$.  
Let $p=(0,0)$ and $q=(2,0)$. Let $e_n=(1,\frac{1}{n})$ for each $n\in\NN$. 
Let $X=\{pq\}\cup\{pe_n\}_{n=1}^\infty\cup\{e_nq\}_{n=1}^\infty$. Note that $X$ is 
a nonlocally connected subcontinuum of $\RR^2$.

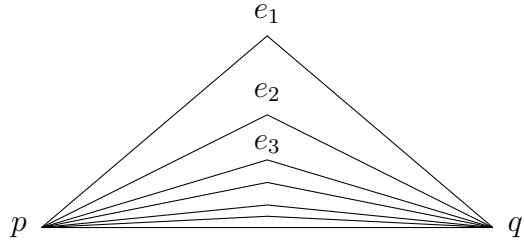
\begin{figure}[h]
    \centering
    \begin{tikzpicture}[scale = 3]
    \draw (0,0) -- (2,0);
    \draw (0,0) -- (1,0.85) -- (2,0);
    \draw (0,0) -- (1,0.5) -- (2,0);
    \draw (0,0) -- (1,0.3) -- (2,0);
    \draw (0,0) -- (1,0.2) -- (2,0);
    \draw (0,0) -- (1,0.1) -- (2,0);
    \draw (0,0) -- (1,0.05) -- (2,0);
    \draw (-0.1,0) node {$p$};
    \draw (2.1,0) node {$q$};
    \draw (1,0.95) node {$e_1$};
    \draw (1,0.6) node {$e_2$};
    \draw (1,0.37) node {$e_3$};
    \end{tikzpicture}
    \caption{Harmonic suspension}\label{fig001}
\end{figure}

We prove that $\FF_1(X)$ is a continuum of colocal connectedness of $\FF_2(X)$. 
Let $\mathcal{U}$ be an open subset of $\FF_2(X)$ such that 
$\FF_1(X)\subseteq \mathcal{U}$. Since $\FF_1(X)$ is compact, there exists 
$\varepsilon>0$ such that $\V_{\varepsilon}(\FF_1(X))\subseteq \mathcal{U}$. 
Let $\{a,b\}\in \FF_2(X)\setminus \V_{\varepsilon}(\FF_1(X))$. Note that 
$\V_{\varepsilon}(a)\cap\V_{\varepsilon}(b)=\emptyset$. Let $A_p$ be the arc 
joining $a$ and $p$ where $q\notin A_p$; let $A_q$ be the arc joining $a$ and $q$ 
where $p\notin A_q$; let $B_p$ be the arc joining $b$ and $p$ where $q\notin B_p$ and 
let $B_q$ be the arc with end points $b$ and $q$ such that $p\notin B_q$.
Observe that $\langle A_p, B_q\rangle\cap\FF_2(X)$ and 
$\langle A_q, B_p\rangle\cap\FF_2(X)$ are subcontinua of $\FF_2(X)$ containing 
$\{\{a,b\},\{p,q\}\}$ \cite[Lemma 1]{JMM}. Also, we have that either 
$\V_{\varepsilon}(A_p)\cap \V_{\varepsilon}(B_q)=\emptyset$ or 
$\V_{\varepsilon}(A_q)\cap \V_{\varepsilon}(B_p)=\emptyset$; i.e., 
$\V_{\varepsilon}(\langle A_q, B_p\rangle)\cap\FF_1(X)=\emptyset$ or 
$\V_{\varepsilon}(\langle A_p, B_q\rangle)\cap\FF_1(X)=\emptyset$. In any case, 
we have that there exists a continuum in 
$\FF_2(X)\setminus \V_{\varepsilon}(\FF_1(X))$ containing the set 
$\{\{a,b\},\{p,q\}\}$. Thus, $\FF_2(X)\setminus \V_{\varepsilon}(\FF_1(X))$ 
is connected. Therefore, $\FF_1(X)$ is a continuum of colocal connectedness of $\FF_2(X)$.
\par\Q

\section{Results on $\T$-closed sets}\label{Tclosed}

In this section we prove general properties of $\T$-closed subsets of continua.
We start with the following characterization theorem.

\begin{theorem}\label{charTclosedset}
Let $X$ be a continuum and let $A$ be a nonempty closed subset of $X$. 
Then the following statements are equivalent:
\smallskip

$(1)$ $A$ is in $\mathfrak{T}(X)$.
\smallskip

$(2)$ For each open subset $U$ of $X$ containing $A$, there exists
an open subset $V$ of $X$ such that $A\subset V\subset U$ and $X\setminus V$
has only finitely many components.
\smallskip

$(3)$ For every subcontinuum $B$ of $X$, disjoint from $A$, 
there exists a subcontinuum
$W$ of $X$ such that $B\subset Int(W)\subset W\subset X\setminus A$.
\end{theorem}

{\bf Proof.} Note that, by \cite[Theorem 4.2.1]{MS}, $(1)$ and $(2)$ 
are equivalent.
Suppose $(2)$ and let $B$ be a subcontinuum of $X$ disjoint from $A$. Let $U$ be
an open subset of $X$ containing $A$ and $Cl(U)\cap B=\emptyset$. 
By our assumption, there exists an open
subset $V$ of $X$ such that $A\subset V\subset U$ and $X\setminus V$ has only
finitely many components. Let $W$ be the component of $X\setminus V$ containing
$B$. Then $W$ is a subcontinuum of $X$ such that 
$B\subset Int(W)\subset W\subset X\setminus A$ (the first inclusion follows from
\cite[Lemma 1.4.1]{MS}). Hence, $(2)$ implies $(3)$.
Next, assume $(3)$ and let $x\in X\setminus A$. Since $\{x\}$ is a 
subcontinuum of $X$,
by our hypothesis, there exists a subcontinuum $W$ of $X$ such that 
$x\in Int(W)\subset W\subset X\setminus A$. Hence, $x\in X\setminus\T(A)$.
Thus, $\T(A)\subset A$ and $A\in\mathfrak{T}(X)$. Hence, $(3)$ 
implies $(1)$.\par\Q

\begin{remark}\label{rem001}
{\rm Regarding \cite[Problem 6.7]{CCOR}, we note that 
\cite[Theorem 4.7]{MKT} gives a partial answer to
this question, since it is shown the following: Let $X$ be a continuum 
with the property of Kelley. Then 
$\mathfrak{T}_\C(X)$ is a subcontinuum of $\C_1(X)$ if and only if 
$\T$ is continuous on continua.
It is an open question if the continuity of $\T$ on continua implies 
the continuity of $\T$ \cite[Question 8.1.10]{MS}. Also,
if the continuum $X$ is the product of two continua, then the answer 
to \cite[Problem 6.7]{CCOR}
is positive since, in this case, $X$ is locally connected 
\cite[Theorem 5.3.6]{MS}.}
\end{remark}

\begin{remark}\label{rem002}
{\rm Regarding \cite[Question 3.23]{OPQV}, we observe that 
\cite[Theorem 4.7]{MKT} gives a partial answer to
this question, since the following is proven: Let $X$ be a continuum 
with the property of Kelley. Then 
$\mathfrak{T}(X)$ is a subcontinuum of $2^X$ if and only if $\T$ is 
continuous. Also,
if the continuum $X$ is the product of two continua, then the answer 
to \cite[Question 3.23]{OPQV}
is positive since, in this case, $X$ is locally connected 
\cite[Theorem 5.3.6]{MS}.}
\end{remark}

\begin{theorem}\label{theo0001}
Let $X$ and $Y$ be chainable continua and let $A$ and $B$ be 
terminal proper subcontinua of $X$ and $Y$, respectively. If 
$C\in\C_1(X\times Y)$ is such that $\pi_X(C)=A$ and 
$\pi_Y(C)=B$, then $\T_{X\times Y}(C)=A\times B$.
\end{theorem}

\textbf{Proof.} Since $C\subseteq A\times B$ and 
$\T_{X\times Y}(A\times B)=A\times B$ \cite[Theorem 4.1.5]{MS}, we have
that $\T_{X\times Y}(C)\subseteq A\times B$. Suppose that there 
exists $(p,q)\in (A\times B)\setminus \T_{X\times Y}(C)$. Then 
there exists a continuum $W$ such that $(p,q)\in Int(W)$ and 
$W\cap C=\emptyset$. Note that $A\subsetneq \pi_X(W)$ and 
$B\subsetneq \pi_Y(W)$, because $Int(W)\neq\emptyset$, both
$\pi_X(W)\cap A\neq\emptyset$ and $\pi_Y(W)\cap 
B\neq\emptyset$, and terminal subcontinua have empty interior
\cite[Lemma 1.4.51]{MS}. 

Let $\alpha\colon [0,1]\to\C_1(W)$ be an order arc such that
$\alpha(0)=\{(p,q)\}$ and $\alpha(1)=W$ \cite[Theorem (1.8)]{N1}. 
Since $A$ and $B$ are
terminal, there exist $t,s\in [0,1]$ such that 
$$\pi_X(\alpha(t))=A \text{ and }\pi_Y(\alpha(s))=B.$$ 
Without
loss of generality we assume that $t\leq s$. Let $D=\alpha(t)$.
Note that $D\subset A\times B$, $\pi_X(D)=A$, 
$D\cap C=\emptyset$ and $\pi_Y(C)=B$. Since $X$ and $Y$ are
chainable 
continua, we contradict Theorem \ref{theo0}. Therefore, 
$\T_{X\times Y}(C)=A\times B$.\par\Q

\begin{theorem}\label{theo0003}
Let $X$ be the pseudo-arc. If $C\in\mathfrak{T}_\C(X^2)$, then either:
\smallskip

$(1)$  $C=X^2$, or
\smallskip

$(2)$ $C=A\times B$, where $A$ and $B$ are proper subcontinua of $X$.
\end{theorem}

\textbf{Proof.} Since $X$ is homogeneous \cite[Theorem 13]{Bing}, 
we have that
$\T_{X^2}$ is idempotent on closed sets
\cite[Theorem 4.2.32]{MT}. Hence, 
$$\mathfrak{T}_\C(X^2)=\{\T_{X^2}(C)\ |\ C\in\C_1(X^2)\}.$$
Also, $X$ is chainable \cite[Definition 7]{Moise} and 
every proper subcontinuum is 
terminal. Thus, the theorem follows from 
Theorem~\ref{theo0001}.\par\Q

\begin{theorem}\label{AnotTclosedAxYnotTclosed}
Let $X$ be a continuum and let $A$ be a nonempty closed subset of $X$
such that $\T_X(A)\not=A$. If $Y$ is a continuum, then $A\times Y$ is not
a $\T$-closed set in $X\times Y$. 
\end{theorem}

{\bf Proof.} Assume that $A\times Y$ is a $\T$-closed set in $X\times Y$.
Let $x\in\T(A)\setminus A$ and $y\in Y$. Then there exists a subcontinuum
$W$ of $X\times Y$ such that $(x,y)\in Int_{X\times Y}(W)\subset W\subset 
(X\times Y)\setminus (A\times Y)$. Note that $\pi_1(W)$ is a subcontinuum
of $X$ and $x\in Int_X(\pi_1(W))$. Hence, $\pi(W)\cap A\not=\emptyset$. 
Thus, there exists $(z,w)\in W$ such that $z\in A$. This implies that
$(z,w)\in W\cap (A\times Y)$, a contradiction. Therefore, $A\times Y$ is
not a $\T$-closed set in $X\times Y$.\par\Q

\begin{theorem}\label{TclosedcompactiffXYlocconn}
Let $X$ and $Y$ be continua. Then $\mathfrak{T}(X\times Y)$ is compact if and
only if $X\times Y$ is locally connected.
\end{theorem}

{\bf Proof.} If $X\times Y$ is locally connected, by \cite[Theorem 2.1.37]{MS},
$\mathfrak{T}(X\times Y)=2^X$. Hence, it is compact \cite[Theorem 1.8.5]{MT}.

Suppose that $X\times Y$ is not locally connected. Without loss of generality,
we assume that $X$ is not locally connected. Thus, there exists a nonempty
closed subset $A$ of $X$ such that $\T(A)\not=A$. Hence, 
$A\times Y$ is not a $\T$-closed set in $X\times Y$ 
(Theorem~\ref{AnotTclosedAxYnotTclosed}). Now, let $\{L_n\}_{n=1}^\infty$
be a sequence of proper subcontinua of $Y$ converging to $Y$. Observe that,
on one hand, $\T_{X\times Y}(A\times L_n)=A\times L_n$, for all $n\in\NN$ 
\cite[Theorem 4.1.5]{MS}. On the other hand, the sequence
$\{A\times L_n\}_{n=1}^\infty$ converges to $A\times Y$. Therefore,
$\mathfrak{T}(X\times Y)$ is not compact.\par\Q

\begin{corollary}\label{TclosedcompactiffXYlocconn2}
Let $X$ and $Y$ be continua. Then $\mathfrak{T}_\C(X\times Y)$ is compact if and
only if $X\times Y$ is locally connected.
\end{corollary}

{\bf Proof.} If $X\times Y$ is locally connected, by \cite[Theorem 2.1.38]{MS},
$\mathfrak{T}_\C(X\times Y)=\C_1(X)$. Hence, it is compact 
\cite[Theorem 1.8.5]{MT}.
The proof of the converse implication is similar to the one given in 
Theorem~\ref{TclosedcompactiffXYlocconn}.\par\Q

\begin{remark}
{\rm Note that in the first line of the proof of \cite[Corollary 4.6]{CCOR}, 
the authors claim
that $\mathfrak{T}_\C(X\times [0,1])$ is a continuum. By 
Corollary~\ref{TclosedcompactiffXYlocconn2}, this is only true when $X$ is
a locally connected continuum. However, the proof of \cite[Corollary 4.6]{CCOR} is correct if the sentence is
changed to $\mathfrak{T}_\C(X\times [0,1])$ is connected.}
\end{remark}

The next corollary follows from Theorem~\ref{TclosedcompactiffXYlocconn} and
Corollary~\ref{TclosedcompactiffXYlocconn2}. It gives a partial answer to 
\cite[Problem 6.7]{CCOR} and \cite[Question 3.23]{OPQV}.

\begin{corollary}
Let $X$ and $Y$ be continua. Then
\smallskip

$(1)$ $\mathfrak{T}_\C(X\times Y)$ is compact if and only if $\T_{X\times Y}$ 
is continuous. 
\smallskip

$(2)$ $\mathfrak{T}(X\times Y)$ is compact if and only if $\T_{X\times Y}$ 
is continuous. 
\end{corollary}

\begin{theorem}\label{pseudoarcTclosedconn}
If $X$ is the pseudo-arc, then $\mathfrak{T}_\C(X^2)$ 
is arcwise connected and not compact.
\end{theorem}

\textbf{Proof.} We see that $\mathfrak{T}_\C(X^2)$ is arcwise 
connected. Let $C\in\mathfrak{T}_\C(X^2)\setminus \{X^2\}$. By 
Theorem~\ref{theo0003}, $C=A\times B$, where $A$ and $B$ are 
proper subcontinua of $X$. Let 
$\alpha,\lambda\colon [0,1] \to\C_1(X)$ be order arcs such that 
$\alpha(0)=A$, $\lambda(0)=B$ and 
$\alpha(1)=\lambda(1)=X$ \cite[Theorem (1.8)]{N1}. Let 
$\beta\colon [0,1]\to\mathfrak{T}_\C(X^2)$ 
be given by $\beta(t)=\alpha(t)\times \lambda(t)$, for each
$t\in [0,1]$. Note that $\beta$ is a map. Therefore, 
$\mathfrak{T}_\C(X^2)$ is arcwise connected.
The fact that $\mathfrak{T}_\C(X^2)$ is not compact, follows from
Theorem~\ref{TclosedcompactiffXYlocconn2}.\par\Q


The following example shows that Theorem \ref{theo0001} does 
not hold, if $X$ is irreducible, but it is not chainable.

\begin{example}
{\rm Let $X=S\cup R$, where $R$ is homeomorphic to $[0,1)$, 
$Cl(R)\setminus R=S$ and $S=\{z\in\mathbb C\ |\ 
|z|=1\}$, as we show in Figure \ref{fig001} (this space is also considered in 
\cite[Example 32]{IMMM2}). Observe that $S$ is a proper
terminal subcontinuum of $X$.
Let $\Delta_S=\{(z,z)\in X^2\ |\ z\in S\}$. Observe 
that $\pi_1(\Delta_S)=\pi_2(\Delta_S)=S$. 
Also, by \cite[Example 32]{IMMM2}, for each 
$(x,y)\in S^2\setminus \Delta_S$ there 
exists an open subset $U$ of $X^2$ such that $\Delta_S\subseteq U$, 
$(x,y)\in X^2\setminus Cl(U)$ 
and $X^2\setminus U$ is a continuum. Therefore, $\T(\Delta_S)=\Delta_S$ 
and $\T(\Delta_S)\neq S^2$. Note that, in this case, $\Delta_S\in\mathfrak{T}_\C(X^2)$.}

\begin{figure}[h]
    \centering
    \begin{tikzpicture}[scale = 1.2]
    \draw (0,0) circle (1cm);
    \draw (1.8,0) arc (0:180:1.65);
    \draw (-1.5,0) arc (-180:0:1.4);
    \draw (1.3,0) arc (0:180:1.3);
    \draw (-1.3,0) arc (-180:0:1.2);
    \draw (1.1,0) arc (0:175:1.08);
    \draw (-1.05,0) node {\small{$\vdots$}};
    \draw (0.6,0) node {$S$};
    \draw (2,1) node {$R$};
    
    \end{tikzpicture}
    \caption{ Spiral}\label{fig001}
\end{figure}
\end{example}

\begin{theorem}
Let $X$ and $Y$ be continua and let $A$ and $B$ be proper, nondegenerate, 
terminal subcontinua of $X$ and $Y$,  respectively. If $a\in A$ and 
$b\in B$, then 
$\T((\{a\}\times B)\cup(A\times\{b\}))=A\times B$.
\end{theorem}

{\bf Proof.} We know that $\T(A\times B)=A\times B$ 
\cite[Theorem 4.1.5]{MS}. Since $(\{a\}\times B)\cup(A\times\{b\})\subseteq A\times B$, 
$\T((\{a\}\times B)\cup(A\times\{b\}))\subset \T(A\times B)$. Hence,
$\T((\{a\}\times B)\cup(A\times\{b\}))\subset A\times B$. 
Suppose that there exists
$(x,y)\in A\times B\setminus\T((\{a\}\times B)\cup(A\times\{b\}))$. Then there
exists a subcontinuum $W$ of $X\times Y$ such that 
$(x,y)\in Int_{X\times Y}(W)\subset W\subset [(X\times Y)\setminus 
(\{a\}\times B)\cup(A\times\{b\})]$.
To obtain a contradiction, we prove that $W\cap (\{a\}\times B)\cup(A\times
\{b\})\not=\emptyset$.
To this end, let $\alpha\colon [0,1]\to\C_1(X\times Y)$ be an order arc such that
$\alpha(0)=\{(x,y)\}$ and $\alpha(1)=W$ \cite[Theorem (1.8)]{N1}. Let 
$\pi_1\colon X\times Y\tto X$ and $\pi_2\colon X\times Y\tto Y$ 
be the projection maps.
Note that  $\pi_1(\alpha(0))=\{x\}\subset A$ and 
$\pi_2(\alpha(0))=\{y\}\subset B$.
Let
$$t_1=\max\{t\in [0,1]\ |\ \pi_1(\alpha(t))\subset A\}$$
and
$$t_2=\max\{t\in [0,1]\ |\ \pi_2(\alpha(t))\subset B\}.$$
Since $\pi_1(W)\cap A\not=\emptyset$, $Int_X(A)=\emptyset$ 
\cite[Lemma 1.4.51]{MS}
($A$ is a terminal subcontinuum of $X$) and $Int_X(\pi_1(W)\not=\emptyset$, 
we have that $A\subset\pi_1(W)$. Similarly,
$B\subset\pi_2(W)$. Assume that $t_1\leq t_2$. Hence, $\pi_1(\alpha(t_1))=A$ and
$\pi_2(\alpha(t_1))\subset B$. From here, we have that $a\in\pi(\alpha(t_1))$.
Thus, there
exists $y'\in\pi_2(\alpha(t_1))$ such that $(a,y')\in\alpha(t_1)\subset W$. 
Since $\pi_2(\alpha(t_1))\subset B$, we obtain that $(a,y')\in\{a\}\times B$.
Hence, $W\cap (\{a\}\times B)\not=\emptyset$, a contradiction. We have a
similar situation if $t_2\leq t_1$. Therefore, 
$\T((\{a\}\times B)\cup(A\times\{b\}))=A\times B$.\par\Q

\begin{corollary}
Let $X$ and $Y$ be continua and let $A$ and $B$ be terminal subcontinua 
of $X$ and $Y$, 
respectively. If $C$ is a subcontinuum of $A\times B$ such that 
$(\{a\}\times B)\cup(A\times\{b\})\subseteq C$, for some $a\in A$ and $b\in B$, 
then $\T(C)=A\times B$. 
\end{corollary}

\textbf{Proof.} Since $(\{a\}\times B)\cup(A\times\{b\})\subseteq C\subseteq A\times B$ 
and $\T((\{a\}\times B)\cup(A\times\{b\}))=\T(A\times B)=A\times B$, 
we have that $\T(C)=A\times B$. \par\Q

\begin{theorem}\label{ATclosedequiv}
Let $X$ be a continuum and let $A$ be a $\T$-closed subcontinuum of $X$
such that $Int_X(A)=\emptyset$. Then 
the following are equivalent:
\smallskip

$(1)$ $A$ is a subcontinuum of colocal connectedness of $X$.
\smallskip

$(2)$ $X\setminus A$ is strongly continuumwise connected.
\smallskip

$(3)$ $X\setminus A$ is continuumwise connected.
\smallskip

$(4)$ $A$ is a nonblock subcontinuum of $X$.
\smallskip

$(5)$ $A$ is a shore subcontinuum of $X$.
\smallskip

$(6)$ $A$ is not a strong centre in $X$.
\end{theorem}

{\bf Proof.} By Lemma~\ref{lemyhfgt5}, we only need to show that (6) implies (1). Let $U$ be an open subset of $X$ containing $A$. 
Since $A$ is a $\T$-closed subset of $X$, by Theorem~\ref{charTclosedset}, there
exists an open subset $V'$ of $X$ such that $A\subset V'\subset U$ 
and $X\setminus V'$
has only finitely many components. Let $N_1,\ldots, N_k$ be the components of 
$X\setminus V'$. By Lemma~\ref{lem7yg6}, without loss of generality,
we assume that $Int_X(N_j)\not=\emptyset$. Since $A$ is not a strong 
centre, for every
$j\in\{1,\ldots, k-1\}$, there exists a subcontinuum $W_j$ of $X\setminus A$ such
that $W_j\cap Int_X(N_j)\not=\emptyset$ and 
$W_j\cap Int_X(N_{j+1})\not=\emptyset$. Let
$W=\left(\bigcup_{j=1}^kN_j\right)\cup\left(\bigcup_{j=1}^{k-1}W_j\right)$, and 
let $V=X\setminus W$. Then $V$ is an open subset of $X$ $A\subset V\subset U$ and
$X\setminus V$ is connected. Therefore,  $A$ is a subcontinuum of 
colocal connectedness of $X$.\par\Q

\begin{remark}
{\rm Let$X$ be a continuum. In Corollary~\ref{scimpliesnonblock2}, we present the equivalence between
$(4)$, $(5)$ and $(6)$ in Theorem~\ref{ATclosedequiv}, when 
$A$ is the diagonal, $\Delta_X$, of $X^2$,
without assuming that $\Delta_X$ is a $\T$-closed subcontinuum of $X^2$.}    
\end{remark}
 
\section{$\Delta_X$ is a $\T$-closed.}\label{DeltaTclosed}

In this section we restrict our study to properties when 
$\Delta_X$ is or is not a $\T$-closed subcontinuum of $X^2$.

\begin{remark}\label{rem003}
If $X$ is a locally connected continuum, then $\Delta_X$ is a 
$\T$-closed subset of $X^2$ \cite[Theorem 2.1.37]{MS}.
\end{remark}

As a consequence of Theorem~\ref{theodescdiagnequ}, we have:

\begin{theorem}\label{theodescdiagnequ2}
If $X$ is a decomposable continuum, then $\T(\Delta_X)\not=X^2$.
\end{theorem}

\begin{theorem}\label{ctheochainterminal}
Let $X$ be a chainable continuum. If there exists a nondegenerate terminal 
subcontinuum $K$ of $X$, then $\Delta_X$ is not a $\T$-closed subcontinuum
of $X^2$. 
\end{theorem}

{\bf Proof.} By Theorem~\ref{theo0001}, $\T(\Delta_K)=K\times K$. 
Since $\Delta_K\subseteq \Delta_X$, $K\times K\subseteq \T(\Delta_X)$. 
Since $K$ is nondegenerate, $\T(\Delta_X)\neq \Delta_X$. 
Therefore, $\Delta_X$ is not a $\T$-closed subset of $X^2$.\par\Q

\begin{corollary}
Let $X$ and $Y$ be continua, where $Y$ is locally connected, and let
$f\colon X\tto Y$ be a monotone open map. If $Y$ is not an arc and
$\Delta_X$ is a $\T$-closed subset of $X^2$, then
$\Delta_X$ is a subcontinuum of colocal connectedness of $X^2$.
\end{corollary}

{\bf Proof.} By Theorem~\ref{nonblockset}, $\Delta_X$ is nonblock subset of $X^2$.
Then, by
Theorem~\ref{ATclosedequiv}, $\Delta_X$ is a subcontinuum of colocal
connectedness of $X^2$.\par\Q

\begin{theorem}\label{theochaind}
Let $X$ be a chainable continuum. If $X$ is indecomposable, 
then $\T(\Delta_X)=X^2$. Hence, 
$\Delta_X$ is not a $\T$-closed subset of $X^2$.
\end{theorem}

\textbf{Proof.} Note that if $W$ is a subcontinuum of $X^2$ 
such that $Int(W)\neq\emptyset$, 
then both $\pi_1(W)$ and $\pi_2(W)$ are subcontinua of $X$ with 
nonempty interior. 
Since $X$ is indecomposable, $\pi_1(W)=\pi_2(W))=X$ \cite[Corollary 1.4.35]{MS}. 
Thus, $W\cap \Delta_X\neq\emptyset$ (Theorem \ref{theo0}). 
Therefore, $\T(\Delta_X)=X^2$.\par\Q

Next result follows from Theorems \ref{theodescdiagnequ2} and \ref{theochaind}.

\begin{corollary}\label{XchainableDeltaX2}
Let $X$ be a chainable continuum. Then $\T(\Delta_X)=X^2$ if and only if $X$ is indecomposable.
\end{corollary}

The following theorem tells us that the arc is the only chainable continuum
for which $\Delta_X$ is $\T$-closed subcontinuum of $X^2$.

\begin{theorem}\label{XchainDeltaTclosedXarc}
Let $X$ be a chainable continuum. Then $\Delta_X$ is $\T$-closed subcontinuum
of $X^2$ if and only if $X$ is an arc.
\end{theorem}

{\bf Proof.} Suppose that $\Delta_X$ is $\T$-closed. Since $X$ is chainable, 
there exist two points $p$ and $q$ in $X$ such that $X$ is irreducible 
between $p$ and $q$ 
\cite[Theorem 12.5]{N2}. Note that $X=A\cup B$, where $A$ and $B$ are proper 
subcontinua of $X$ (Theorem \ref{theochaind}). Without loss of 
generality, we may assume that $p\in A\setminus B$ and $q\in B\setminus A$. 
We consider two cases: 
\medskip

\textbf{Case (1).} $\{p\}=X\setminus \kappa(q)$ and 
$\{q\}=X\setminus \kappa(p)$. 
\smallskip

Let
$$M=\{(x,y)\in X^2\ |\ \text{ there exists }K\in\C_1(X^2\setminus\Delta_X)
\text{ such that }\{(p,q),(x,y)\}\subseteq K\}$$ 
and
$$N=\{(x,y)\in X^2\ |\ \text{ there exists }K\in\C_1(X^2\setminus\Delta_X)
\text{ such that }\{(q,p),(x,y)\}\subseteq K\}.$$
By \cite[Lemma 2.1]{IMMM1}, $M\cap N=\emptyset$. Also, $M$ and $N$ 
are continuumwise connected.

We prove that $Cl(M\cup N)=X^2$. Let $\alpha,\beta\colon [0,1]\to\C_1(X)$ be 
order arcs such that $\alpha(0)=\{p\}$, $\beta(0)=\{q\}$ and 
$\alpha(1)=\beta(1)=X$
\cite[Theorem (1.8)]{N1}. 
Let $(z,w)\in X^2$ and let $U$ and $V$ be open subsets of $X$ such 
that $(z,w)\in U\times V$. 
We consider two more subcases:
\medskip

\textbf{Subcase (1.1).} There exist $t<s$ such that 
$\alpha(t)\cap U\neq\emptyset$ and 
$(\alpha(s)\setminus\alpha(t))\cap V\neq\emptyset$. 
\smallskip

Let $a\in \alpha(t)\cap U$ and let 
$b\in (\alpha(s)\setminus\alpha(t))\cap V$. Let $I=\alpha(t)\times \{b\}$. 
Observe that 
$(p,b)\in I$ and $I\in\C_1(X^2\setminus\Delta_X)$. Since 
$\{p\}=X\setminus \kappa(q)$, 
there exists a proper subcontinuum $Z$ of $X$ such that $\{b,q\}\subseteq Z$.
Note 
that $p\in X\setminus Z$ and $J=\{p\}\times Z$ is a subcontinuum 
of $X^2\setminus\Delta_X$. 
Since $(p,b)\in I\cap J$, $I\cup J$ is a subcontinuum of 
$X^2\setminus\Delta_X$ such that 
$(p,q)\in I\cup J$. Thus, $I\cup J\subseteq M$. Since 
$(I\cup J)\cap (U\times V)\neq\emptyset$, 
we have that $(z,w)\in Cl(M)$.
\medskip

\textbf{Subcase (1.2).} $V\subseteq \alpha(t)$, for each $t$ 
such that $\alpha(t)\cap U\neq\emptyset$. 
\smallskip

Note that there exist $t'<s'$ such that $\beta(t')\cap U\neq\emptyset$ and 
$(\beta(s')\setminus\beta(t'))\cap V\neq\emptyset$; because, if 
$V\subseteq\beta(t')$, 
for each $t$ such that $\beta(t')\cap U\neq\emptyset$, then 
$\alpha(t)\cup\beta(t')=X$ 
and $U\nsubseteq \alpha(t)\cup\beta(t')$, a contradiction. Now, 
with a similar argument 
to the one given in Subcase $(1.1)$, we have that $(z,w)\in Cl(N)$.

Therefore, we have that $Cl(M\cup N)=X^2$. Now, suppose that there 
exists $(r,s)\in Cl(M)\cap N$. 
Since $\Delta_X$ is $\T$-closed, $(r,s)\in N$ and $N\cap\Delta_X=\emptyset$, we 
have that there exists a subcontinuum $W$ such that $(r,s)\in Int(W)$ and 
$W\cap \Delta_X=\emptyset$. Observe that $W\cap M\neq\emptyset$ and 
$W\cap N\neq\emptyset$. Thus, we have both $W\subseteq M$ and $W\subseteq N$,
a contradiction to the fact that $M\cap N=\emptyset$. Therefore, 
$Cl(M)\cap N=\emptyset$. 
Similarly, we have that $M\cap Cl(N)=\emptyset$. Hence, 
$X^2=M\cup N\cup \Delta_X$, where $M$ and $N$ are disjoint open 
subsets of $X^2$. 
By \cite[Theorem 3]{Kat}, $X$ is an arc.

\medskip

\textbf{Case (2).} Either $X\setminus\kappa(p)$ is nondegenerate 
or $X\setminus\kappa(q)$ is nondegenerate.
\smallskip

Assume that $X\setminus\kappa(p)$ is nondegenerate. Let 
$K=Cl(X\setminus\kappa(p))$. 
Observe that $K$ is a nondegenerate proper subcontinuum of $X$. Note that if $K$ 
has empty interior, by Lemma~\ref{lemmat}, $K$ is a terminal 
subcontinuum of $X$. Hence, 
by Corollary~\ref{ctheochainterminal}, $\Delta_X$ is not 
$\T$-closed. Thus, suppose that $Int(K)\neq\emptyset$. 
By Lemma \ref{lemmat}, $K$ is indecomposable. Let 
$(a,b)\in X^2\setminus\Delta_X$ be such that 
$a$ and $b$ are in $K$. Since $\Delta_X$ is $\T$-closed, there 
exists a subcontinuum $W$ of $X^2$ such 
that $(a,b)\in Int(W)\subseteq W\subseteq X^2\setminus\Delta_X$. Note that 
$K\subseteq \pi_1(W)\cap \pi_2(W)$ and, $\pi_1(W)\subseteq \pi_2(W)$ or 
$\pi_2(W)\subseteq \pi_1(W)$. Suppose that $\pi_1(W)\subseteq\pi_2(W)$ and let 
$L=\pi_2(W)$. Since $X$ is chainable, $L$ is chainable \cite[Theorem (9.C.4)]{Chris}. 
Moreover, $\Delta_L$ 
and $W$ are subcontinua of $L^2$ such that $\pi_1(\Delta_L)=L$ and $\pi_2(W)=L$. 
Thus, by Theorem \ref{theo0}, $W\cap \Delta_L\neq\emptyset$. Since 
$\Delta_L\subseteq\Delta_X$, $W\cap \Delta_X\neq\emptyset$. A contradiction. 
Therefore, by Cases $(1)$ and $(2)$, we have that $X$ is an arc.

Conversely, if $X$ is an arc, then $X^2$ is locally connected. Hence, by 
\cite[Theorem 2.1.37]{MS}, $\Delta_X$ is $\T$-closed.\par\Q

\begin{theorem}\label{Deltachar}
Let $X$ be a continuum that is not an arc, and consider the following;
\smallskip

$(1)$  $\Delta_X$ is a $\T$-closed subcontinuum of $X^2$.
\smallskip

$(2)$ $X^2\setminus\Delta_X$ is strongly continuumwise connected.
\smallskip

$(3)$ $\Delta_X$ is a subcontinuum of colocal connectedness of $X^2$.
\smallskip

$(4)$ $X^2\setminus\Delta_X$ is continuumwise connected.
\smallskip

\noindent Then $(1)$, $(2)$ and $(3)$ are equivalent, $(1)$ implies $(4)$. 
Moreover, if $X$ has the property of Kelley, then $(4)$ is equivalent $(1)$.
\end{theorem}

{\bf Proof.} Suppose $\Delta_X$ is a $\T$-closed subcontinuum 
of $X^2$. Then, since $X$ is not an arc, by
\cite[Theorem 3]{Kat}, $X^2\setminus\Delta_X$ is connected.
Hence, by
\cite[Theorem 4.2.9]{MS}, $X^2\setminus\Delta_X$ is an 
exhausted $\sigma$-continuum.
Now, by \cite[Theorem 1.4.75 $(2)$]{MS}, $X^2\setminus\Delta_X$ 
is strongly continuumwise connected.
Assume that $X^2\setminus\Delta_X$ is strongly continuumwise 
connected.
Since $X^2\setminus\Delta_X$ is open, we have that 
$X^2\setminus\Delta_X$
is open and strongly continuumwise connected. By 
\cite[Theorem 1.4.75 $(2)$]{MS},
$X^2\setminus\Delta_X$ is an exhausted $\sigma$-continuum. 
Therefore,
by \cite[Theorem 4.2.9]{MS}, $\Delta_X$ is a $\T$-closed 
subcontinuum of $X^2$. Hence, $(1)$ and $(2)$ are equivalent.

By Theorem \ref{lemma8uhgvfr5}, we have that $(2)$ and $(3)$ are equivalent.

If $\Delta_X$ is a $\T$-closed subcontinuum of $X^2$, by 
Theorem~\ref{ATclosedequiv}, $X^2\setminus\Delta_X$ is continuumwise connected.
Thus, $(1)$ implies $(4)$. Suppose $X$ has the property of Kelley and
$X^2\setminus\Delta_X$ is continuumwise connected. Since $X$ has the
property of Kelley, by Theorem~\ref{strongcontconniffcontconn},
$X^2\setminus\Delta_X$ is strongly continuumwise connected.
Therefore,  
$\Delta_X$ is a $\T$-closed subcontinuum of $X^2$.\par\Q

\begin{remark}
By \cite[Corollary 15]{IMMM2}, if $X$ is the double harmonic fan, 
then $X^2\setminus\Delta_X$ is continuumwise connected, but 
$\Delta_X$ is not a subcontinuum of colocal connectedness of $X^2$. 
Hence, by Theorem~\ref{Deltachar}, $\Delta_X$ is not $\T$-closed. 
This shows that property if Kelley is needed in Theorem~\ref{Deltachar}.
\end{remark}

\begin{remark}
Note that, by Theorem~\ref{XchainDeltaTclosedXarc}, if $X$ is a chainable 
continuum different than an arc, then $\Delta_X$ is not $\T$-closed
subset of $X^2$. Also, 
by \cite[Remark 34]{IMMM2} and Theorem~\ref{Deltachar}, we have that 
$X$ is either a solenoid or the Warsaw circle, then $\Delta_X$ is not $\T$-closed
subcontinuum of $X^2$.  
\end{remark}








\section{Models}\label{Models}

In this section, we present partial answers to \cite[Problem 6.6]{CCOR}.
Note that, since atomic maps are monotone \cite[Theorem 8.1.24]{MT}, 
we can extend \cite[Theorem 6.3]{CCOR} to continua; i.e.,

\begin{theorem}\label{fatomopencont}
Let $X$ and $Y$ be continua and let $f\colon X\tto Y$ be an atomic open map.
Then $\mathfrak{T}_\C(X)$ is homeomorphic to $\mathfrak{T}_\C(Y)$.
\end{theorem}

\begin{theorem}\label{theo87y6}
Let $X$ and $Y$ be continua and let $f\colon X\tto Y$ be a monotone open map.
Then:
\smallskip

$(1)$ $\mathfrak{T}(Y)$ may be embedded in $\mathfrak{T}(X)$
\smallskip

\noindent and
\smallskip

$(2)$ $\mathfrak{T}_\C(Y)$ may be embedded in 
$\mathfrak{T}_\C(X)$.
\end{theorem}

{\bf Proof.} Suppose $f$ is an open map. Hence, by 
\cite[Theorem 1.6.16]{MS}, 
$\Im(f)\colon 2^Y\tto 2^X$ given by $\Im(f)(B)=f^{-1}(B)$ 
is a map. Hence, since
$\Im(f)$ is a one-to-one map between compacta, $\Im(f)$ 
is an embedding.
Since $f$ is monotone too, by \cite[Theorem 4.1.25]{MS}, 
we have that 
$\Im(f)(\mathfrak{T}(Y))\subset\mathfrak{T}(X)$, and
$\Im(f)(\mathfrak{T}_\C(Y))\subset\mathfrak{T}_\C(X)$. Thus, 
$\mathfrak{T}(Y)$ may be embedded in $\mathfrak{T}(X)$, and
$\mathfrak{T}_\C(Y)$ may be embedded in 
$\mathfrak{T}_\C(X)$\par\Q

\begin{remark}
Note that, by \cite[Theorem 4.12]{BFM}, if $X$ is a type $\lambda$ continuum,
then there exists a bijection between $\mathfrak{T}(X)$ and $2^{[0,1]}$, whose
restriction to $\mathfrak{T}_\C(X)$ gives a bijection between 
$\mathfrak{T}_\C(X)$
and $\C_1([0,1])$.
\end{remark}

\begin{theorem}\label{models}
Suppose that a continuum $X$ belongs to the following classes of continua:
\smallskip

$(1)$ Continuously irreducible continua.
\smallskip

$(2)$ Continuously type $A'$ $\theta$-continuum.
\smallskip

$(3)$ Continuously type $A'$ $\theta_n$-continuum, $n\in\NN$.
\smallskip

\noindent Then we have that:
\smallskip

$(1')$ $\mathfrak{T}(X)$ is homeomorphic to the Hilbert cube and
$\mathfrak{T}_\C(X)$ is homeomorphic to a $2$-cell.
\smallskip

$(2')$ $\mathfrak{T}(X)$ is homeomorphic to the Hilbert cube and
$\mathfrak{T}_\C(X)$ is homeomorphic to a polyhedron.
\smallskip

$(3')$ $\mathfrak{T}(X)$ is homeomorphic to the Hilbert cube and
$\mathfrak{T}_\C(X)$ is homeomorphic to a polyhedron.
\end{theorem}

{\bf Proof.} Suppose $(1)$, we show that $(1')$ is true. Since $X$
is a a continuously irreducible continuum, by \cite[Lemma 3.3]{M3},
the quotient map $q\colon X\tto [0,1]$ is an open atomic map.
Then, by \cite[Theorem 6.3]{CCOR}, we have that $\mathfrak{T}(X)$ 
is homeomorphic to $2^{[0,1]}$; i.e., $\mathfrak{T}(X)$ is homeomorphic 
to the Hilbert cube \cite[Theorem 3.6]{SW1}.
By Theorem~\ref{fatomopencont}, we also obtain that 
$\mathfrak{T}_\C(X)$ is homeomorphic to $\C_1([0,1])$; i.e., $\mathfrak{T}_\C(X)$
is homeomorphic to a $2$-cell \cite[1, p. 267]{Duda}.

Assume $(2)$, we prove $(2')$ holds. Since $X$ is 
a continuously type $A'$ $\theta$-continuum, by \cite[Theorem 3.6]{M1},
the quotient map $q\colon X\tto G$, where $G$ is a graph, is an atomic open map.
Then, by \cite[Theorem 6.3]{CCOR}, we have that $\mathfrak{T}(X)$ 
is homeomorphic 
to $2^D$; i.e., $\mathfrak{T}(X)$ is homeomorphic to the Hilbert cube 
\cite[Theorem 4.4]{SW2}. 
By Theorem~\ref{fatomopencont}, we also obtain that 
$\mathfrak{T}_\C(X)$ is homeomorphic to $\C_1(G)$; i.e., $\mathfrak{T}_\C(X)$
is homeomorphic to a polyhedron \cite[6.4, p. 276]{Duda}.

By \cite[Corollary 3.9]{M1}, if $X$ is a continuously type $A'$
$\theta$-continuum, then there exists $n\in\NN$ such that $X$ is a
$\theta_n$-continuum. Hence, the argument of the previous paragraph holds
for continuously type $A'$ $\theta_n$-continua. Thus, assuming $(3)$, we
obtain $(3')$.\par\Q

\begin{theorem}\label{modelshomog}
Suppose $X$ is:
\smallskip

$(1)$ The Menger universal curve and $\widehat X$ is the
Menger universal curve of pseudo-arcs.
\smallskip

$(2)$ The Sierpi\'nski universal plane curve and $\widehat X$
is the Sierpi\'nski universal plane curve of pseudo-arcs.
\smallskip

$(3)$ The circle and $\widehat X$ is the circle of pseudo-arcs.
\smallskip

$(4)$ The arc and $\widehat X$ is the arc of pseudo-arcs.
\smallskip

\noindent Then we obtain that
\smallskip

$(1')$ Both $\mathfrak{T}(\widehat X)$ and $\mathfrak{T}_\C(\widehat X)$ 
are homeomorphic 
to the Hilbert cube.
\smallskip

$(2')$  Both $\mathfrak{T}(\widehat X)$ and $\mathfrak{T}_\C(\widehat X)$ 
are homeomorphic 
to the Hilbert cube.
\smallskip

$(3')$  $\mathfrak{T}(\widehat X)$ is homeomorphic to the Hilbert cube and
$\mathfrak{T}_\C(\widehat X)$ is homeomorphic to a $2$-cell.
\smallskip

$(4')$  $\mathfrak{T}(\widehat X)$ is homeomorphic to the Hilbert cube and
$\mathfrak{T}_\C(\widehat X)$ is homeomorphic to a $2$-cell.
\end{theorem}

{\bf Proof.} If $X$ is the Menger universal curve and $\widehat X$ is the Menger 
universal curve of pseudo-arcs, by \cite[Corollary 6.5]{CCOR}, 
$\mathfrak{T}(\widehat X)$ 
is homeomorphic to $2^X$ and $\mathfrak{T}_\C(\widehat X)$ is 
homeomorphic to $\C_1(X)$.
Thus, $\mathfrak{T}(\widehat X)$ is homeomorphic to the Hilbert cube
\cite[Theorem 3.2]{CS}, and since $X$ does not have frHilbertee arcs, also 
$\mathfrak{T}_\C(\widehat X)$ is homeomorphic to the Hilbert cube
\cite[Theorem 4.1]{CS}.

A similar argument proves the case of the Sierpi\'nski universal plane curve.

Suppose $X$ is the circle and $\widehat X$ is the circle of pseudo-arcs.
Then, by \cite[Corollary 6.5]{CCOR}, $\mathfrak{T}(\widehat X)$ is 
homeomorphic to $2^X$ and $\mathfrak{T}_\C(\widehat X)$ is homeomorphic to 
$\C_1(X)$. In this case, by \cite[Theorem 3.2]{CS}, $\mathfrak{T}(\widehat X)$
is homeomorphic to the Hilbert cube, but $\mathfrak{T}_\C(\widehat X)$
is homeomorphic to a $2$-cell \cite[2, p. 267]{Duda}.

A similar argument shows the case of the arc of pseudo-arcs, we
use \cite[1, p. 267]{Duda}, to conclude that $\mathfrak{T}_\C(X)$ is 
a $2$-cell.\par\Q

Regarding paper \cite{OPQV}, we have the following observations:
\smallskip

$(1)$ In part $(3)$ of Proposition 4.3, the authors imply that if $\T$ is
continuous on continua, then each subcontinuum of $X$ is of $\T$-finite type. 
We do not know if this is true or not. By the word ``similarly'', we think that 
they believe that the following is true:
The continuity of $\T$ on continua implies that $\T$ is idempotent on 
continua. This is still an open question.

$(2)$ Part $(5)$ of Proposition 4.3 is not accurate. If $X$ is a homogeneous
continuum, then $\T$ is idempotent on closed sets \cite[Theorem 4.2.32]{MT}. 
Hence, for each nonempty closed subset $A$ of $X$, we know that $\T^2(A)=\T(A)$. 

$(3)$ Part $(6)$ of Proposition 4.3 is not accurate. Since the authors are only 
considering nonempty closed subsets of a type $\lambda$ continuum $X$, by
the proof of \cite[Lemma 7.5]{BFM}, $\T^3(A)=\T^2(A)$, for all nonempty closed
subsets
$A$ of $X$.

$(4)$ Regarding Corollary 4.8 and Theorem 4.9, it is well known that 
atomic maps are monotone \cite[Theorem 8.1.24]{MT}. Hence, it not
necessary to ask the map on those results to be monotone.

\vskip.5in

{\bf Acknowledgement.} This work was 
supported by UNAM PASPA-DGAPA.
The first named author thanks the support given 
by 
\textit{La Vicerrector\'{\i}a de Investigaci\'on y Extensi\'on 
de la Universidad Industrial 
de Santander y su Programa de Movilidad}.
The second named author thanks the Universidad Industrial de 
Santander, Colombia, 
for the support given during this research.

\noindent (J. Camargo)\\   
Escuela de Matem\'aticas, Facultad de Ciencias,\\ 
Universidad Industrial de Santander, Ciudad Universitaria,\\
Carrera 27 Calle 9, Bucaramanga,\\
Santander, A. A. 678, COLOMBIA.\\
e-mail: jcamargo@saber.uis.edu.co\\

\noindent (S. Mac\'{\i}as)\\
\noindent Current address:\par
\noindent Escuela de Matem\'aticas, Facultad de Ciencias,\\ 
Universidad Industrial de Santander, Ciudad Universitaria,\\
Carrera 27 Calle 9, Bucaramanga,\\
Santander, A. A. 678, COLOMBIA.\\

\noindent Permanent address:\par
\noindent Instituto de Matem\'aticas,\\ 
Universidad Nacional Aut\'onoma de M\'exico,\\
Circuito Exterior, Ciudad Universitaria,\\
CDMX, C. P. 04510. M\'EXICO.\\ 
e-mail: sergiom@matem.unam.mx \\
e-mail: macias@unam.mx

\end{document}